# Energy Storage Peak Shaving Feasibility: Case Studies in Upstate New York

Thomas H. Ortmeyer
Clarkson University
Potsdam, NY 13699

Tuyen Vu
Clarkson University
Potsdam, NY 13699

*Abstract*—**This paper presents the results of a benefit-cost analysis involving the application of battery energy storage systems (BESS) for three of New York State's municipal electric departments (MEDs). New York's municipal electric utilities generally have allocations to receive hydroelectric energy from the New York Power Authority (NYPA). When the NYPA firm hydro allocation is exceeded, the utility must procure energy from other sources, generally at significantly higher prices. Additionally, the municipal utility pays to transmit the energy from the generation source to the utility gateway. There are multiple ways that these MED's can use energy storage technology to reduce their costs.  This paper develops benefit and cost analyses, and concludes that BESS technology is cost-effective in some cases, and will become more cost-effective in the future as battery prices decrease.**

## I. Introduction

This paper provides study results on the costs and benefits of battery energy storage systems (BESS) for three of New York's municipal electric departments. New York's municipal electric utilities generally have contracts to receive hydroelectric energy from the New York Power Authority (NYPA). These hydro allocations provide low-cost energy to the municipal utility. When the NYPA firm hydro allocation is exceeded, the utility must procure energy from other sources, generally at significantly higher prices. Additionally, the municipal utility pays to transmit the energy from the generation source to the utility gateway. As load growth occurs on the utility, electricity charges can rise steeply as the NYPA allocation is exceeded and the wheeling charges increase. The utility can provide value to their customers by managing their energy and power requirements.

At the same time, rapid changes in battery technology are significantly reducing storage costs and improving the feasibility of peak shaving and load shifting applications. This report assesses the potential for BESS systems to provide a cost reduction to the municipal departments, which can be passed on to their customers.

The project identified three methods through which municipal electric departments can benefit from a BESS system:

1) Peak shaving to increase allocation

2) Optimizing market rate purchases

3) Peak shaving in months below allocation

Section II discusses each of these methods in detail.

## II. Cost Saving Options and Study Parameters

### A. Peak Shaving to Increase Allocation

Each MED is allotted a firm hydro allocation from NYPA. NYPA applies an energy charge for energy within the allocation as well as a demand charge. The current rate for firm hydro energy is $0.00492 per kilowatt-hour (kWh). In months when the MED does not exceed its allocation, all energy is purchased at the firm hydro rate.  In months when the MED exceeds its demand allocation, the firm hydro load share is calculated as the ratio of the firm hydro demand divided by the monthly peak demand.  For every hour of the billing period, the electric department receives firm hydroelectric energy equal to the hourly metered load times the firm hydro load share for the month.  Additional energy is purchased at the market rate determined by NYPA and based on the New York Independent System Operator (NYISO) market price.

### B. Peak Shaving to Reduce Demand Charge in months when allocation is not exceeded.

NYPA bills the departments for both demand and energy, and the demand charge is currently $4.07 per kilowatt (kW), based on the peak hourly demand in the given month. This demand charge is for the monthly kW demand up to the firm hydro

The authors acknowledge the support of the New York State Energy Research and Development Authority (NYSERDA), Project 132705, Scott Larsen, Project Manager.

allocation. When demand goes above the allocation, there are no additional demand charges for demand above the allocation.

## C. Peak Shaving to Reduce Demand Charge in months when allocation is not exceeded.

In periods when the power draw is above the adjusted allocation, the cost of energy to the electric department is based on the hourly real time-weighted integrated (TWI) real-time zonal locational-based marginal cost (LBMP) for the New York Independent System Operator (NYISO) zone in which the department is located. The TWI prices are calculated from the real-time market five-minute prices posted by NYISO.

## D. Wheeling Costs

National Grid wheels energy to Tupper Lake and Lake Placid MED's. These wheeling charges are billed on a flat per kWh basis. As a battery installation would not significantly change the total kWh, these wheeling charges are not considered further in the study. Massena MED is fed directly by NYPA. Their cost can be influenced by a battery installation, and their potential reduction in charges in the 2018 test year was estimated in this study.

## E. Energy Prices

The three municipal electric departments in this study are located in northern New York. Lake Placid and Massena are in NYISO's North Zone, and Tupper Lake is in NYISO's Mohawk Valley Zone. In both zones, there has been a downward trend in the market cost of energy in recent years. Figure1 shows the monthly average prices of the NYISO North Zone Real-Time Weighted Hourly costs from 2007–2010 and 2015–2019. These prices are plotted by month over the course of each year. It is clear from Figure 1 that there has been an overall downward trend in these prices in recent years (2015–2019), as compared to the earlier years of 2007–2010. This trend makes forecasting of future benefits of a BESS installation difficult. There is little expectation that energy costs will begin to see significant increases in northern New York until transmission congestion is addressed.

## F. Study Assumptions

- The peak shaving portions of the study were conducted based on the actual demand reduction achieved by the BESS installation. Note that the BESS real power output and the BESS inverter apparent power rating will need to be higher than the corresponding level of demand reduction—perhaps 20-40% higher, due to errors in the load forecasting process, system losses, allowance for VARs delivered to the AC bus, etc.
- The study assumes that the BESS system will have real-time access to the power drawn by the MED, and that the BESS only delivers sufficient power to clamp the power draw at the peak value over the course of each hourly period when the total demand is above the targeted peak level. This will result in a significant reduction in the battery energy required, as compared to BESS units which deliver rated power throughout a discharge period.
- The energy purchased within the NYPA firm hydro allocation is at a base cost of $4.92 per MWh. The energy purchases above the allocation are made at the NYISO hourly time weighted real-time market rate. As this rate fluctuates by the hour, exact cost savings can only be calculated after the fact.
- Note that the municipal departments are billed for additional charges by NYPA. These charges include the NYISO charges and the Clean Energy Standard (CES) charges. The CES charges are the same, regardless if the energy received is at the firm hydro rate or the market rate. The NYISO charges are not broken out by MW or MWh in the billings. These additional charges are not considered further in this analysis.

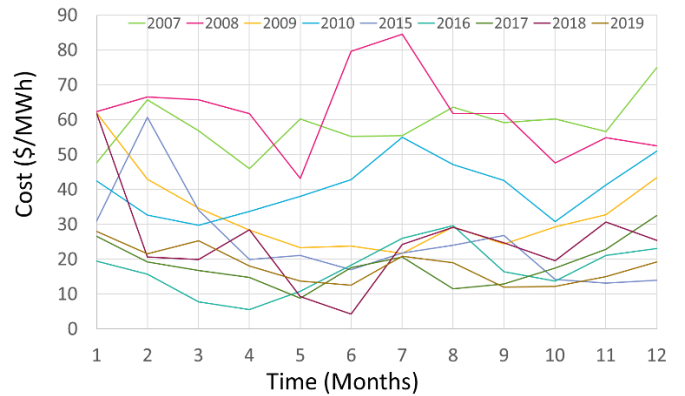

**Figure 1. Monthly Averages for North Zone Hourly Time Weighted Integrated Real-Time Energy Costs by Year.**

## III. ENERGY STORAGE COST STUDY METHOD

The energy storage systems considered in this study will be in the 1MW power range and 1-10 MWh energy range. Therefore, this study focused on the study of energy storage technologies suitable for this size range. The study evaluated 5 battery technologies:
- Lithium nickel manganese cobalt oxide (NMC)
- Lithium iron phosphate (LFP)
- Sodium sulphur (NaS)
- Valve regulated lead acid (VRLA)
- Flooded Lead Acid (FLA)

The BESS system costs were analyzed for the two options identified in the previous section. This analysis is based on BESS cost data from the International Renewable Energy Agency for the reference case **Error! Reference source not found.**.

- The annuity, $C_A$ is calculated as follows:

$$C_A = C_{A\_es} + C_{A\_pc} + C_M + C_{loss}, \qquad (1)$$

where $C_{A\_es}$ is the annualized cost of the energy storage system, $C_{A\_pc}$ is the annualized cost of the power conversion system, $C_M$ is the maintenance cost which is equivalent to 1.5% of annualized costs, and $C_{loss}$ is the loss of electrical-chemical energy conversion depending on the type of battery storage and power conversion module accounting for 2%. The overall loss is converted to *5 cents/kWh* at the day-ahead market rate during the case studies.

The annualized present value, $PV_{annual}$ is calculated as

$$PV_{annual} = \frac{PV \times i}{(1-(1+i)^{-n}) \times (1+i)}, \quad (2)$$

in which $PV$ is the present value of energy storage or power conversion system, $n$ is the number of years that an energy storage system is expected to provide service, and $i$ is the interest rate assuming at 3% annually for the lifetime of the project.

Other factors that can impact the system cost include self-discharge which is small and can be ignored and the depth of discharge (DoD) which depends on the technologies.

IV. CASE I: TUPPER LAKE ELECTRIC MED

Tupper Lake Municipal Electric Department (TLMED) NYPA allocation is 18.845 MW of firm hydro. In 2018, TLMED exceeded its allocation only during January. The yearly peak was 22.6 MW. In 2019, again January was the only month when their allocation was exceeded.

A. Allocation Shifting

Allocation shifting was investigated for both Januarys of the study. BESS installations with 0.5, 1.0, and 2.0 MW of peak shaving capability were investigated. The results are given in Table 1.

**Table 1. Potential for Energy Shifting at Tupper Lake.**

| Peak Shaving capability (MW) | Firm Hydro Purch. (MWh) | Market Purch. (MWh) | Benefit (MWh) | BESS energy Req'd (MWh) |
|---|---|---|---|---|
| January 2018 | | | | |
| 0.0 | 9360 | 1840 | | |
| 0.5 | 9573 | 1628 | 212 | 0.81 |
| 1.0 | 9795 | 1406 | 434 | 6.02 |
| 2.0 | 10271 | 930 | 911 | 26.03 |
| January 2019 | | | | |
| 0.0 | 9798 | 1230 | | |
| 0.5 | 10033 | 1063 | 235 | 0.83 |
| 1.0 | 10280 | 816 | 482 | 2.48 |
| 2.0 | 10811 | 284 | 1013 | 15.40 |

Table 1 shows that the benefits in MWh are consistent between January 2018 and 2019. The benefits increase approximately linearly with peak shaving capability.
The BESS energy required to meet this level of peak shaving, however, increases at a much faster rate. Also, January 2018 was significantly colder than January 2019, as demonstrated by the higher BESS energy required to meet the peak shaving target.

B. Peak Shaving to Reduce Demand Charge

Tupper Lake has the potential to reduce their demand charges in the 11 months when their allocation is not exceeded. Table 2 shows the required BESS energy capability required to meet the peak shaving for the three cases studied.

**Table 2. Required BESS Energy in MWh to Achieve the Targeted Peak Shave in 2018.**

| Month | 0.5 MW peak shave | 1.0 MW peak shave | 2.0 MW peak shave |
|---|---|---|---|
| February | 0.80 | 2.94 | 21.4 |
| March | 0.47 | 1.42 | 4.61 |
| April | 0.57 | 1.82 | 8.93 |
| May | 0.46 | 1.41 | 5.06 |
| June | 1.37 | 6.02 | 39.71 |
| July | 4.21 | 10.93 | 88.71 |
| August | 1.94 | 7.71 | 55.66 |
| September | 0.81 | 4.41 | 18.40 |
| October | 1.26 | 6.11 | 26.02 |
| November | 0.78 | 2.59 | 8.16 |
| December | 0.78 | 2.54 | 9.19 |

Table 2 shows that the energy requirements change significantly between months. For the 0.5 peak shave case, the value ranges from 0.46 MWh to 4.21 MWh. There were four months when the 0.83 MWh was required in January 2019.

In this study, the decision was made to size the BESS energy rating to meet the allocation shifting requirement. In those months when the energy required for peak shaving was greater than this, a lower value of peak shaving would be targeted in order to stay within the energy rating of the installation. Due to the high energy requirements of the 2MW system, it is not considered further.

C. Optimizing Market Rate Purchase

TLMED will have the opportunity to optimize market rate purchases in January when it exceeds the NYPA firm hydro allocation. The service could be conducted on days when the monthly peak would not be established. A study of the 2018 North Zone Day-Ahead market prices showed an average hi/lo price difference of around $25 per MWh, with a range from $0 to $100 per MWh. For this study, a conservative annual savings estimate of $25 per day, per MW of BESS power rating is used. With the BESS unit operating at rated power for one hour of charge and one hour of discharge in 28 days in January, would result in an annual savings of $700. This would not have a significant impact on the BESS unit economics.

### D. Benefit Summary

The total benefits of adding a BESS at Tupper Lake are given in Table 3.

**Table 3. Projected Savings for TLMED in 2018.**

| BESS size | Savings from Allocation Shifting | Savings on Demand Charge | Total Projected Savings |
|---|---|---|---|
| 0.5MW 0.83MWh | $5300 | $19577 | $24,877 |
| 1.0MW 6.0MWh | $10,850 | $42,734 | $53,584 |

### E. Benefits/Costs Analysis

The projected annualized costs are shown in Figure 2, for the 0.5MW, 0.83 MWh BESS. The installation and operating and maintenance costs are determined from Ref. [2], as is the projected system life. The costs are presented as a function of the year installed, using projected cost reductions and performance improvements.

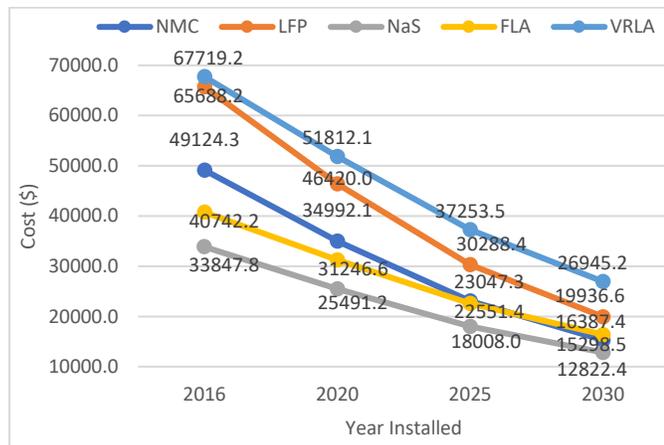

**Figure 2. Projected Annualized Cost for the 0.5MW, 0.83MWh Installation. Costs are based on year of project installation.**

By comparing Table 3 benefits with Figure 2 costs, it can be seen that this installation is near the breakeven point for NaS at present. As prices continue to decline, the economics will improve.

The 1.0MW, 6.0MWh unit, however, is not economically attractive. The cost projection for this case is above the breakeven point through 2030.

## V. CASE 2: LAKE PLACID VILLAGE ELECTRIC DEPARTMENT

The Lake Placid Village Electric Department (LPVED) NYPA allocation is 28.915 MW of firm hydro (28.934 in early 2018). In 2018, LPVED exceeded its allocation for 6 months: January, February, March, April, November, and December. The yearly peak was 51.3 MW, which occurred in January.

### A. Allocation Shifting

Table 4 shows the allocation shifting results for LPVED in 2018 and 2019. The results are largely consistent for these two years. The exception is the energy requirement for the 2.0MW peak shaving case, where the 2018 requirement to meet 2.0MW of peak shaving was 50% higher than for 2019.

In the case of Tupper Lake, the energy requirement increases at a higher rate than the increase in firm hydro purchases as the peak shaving capability of the BESS increases.

**Table 4. Increase in Firm Hydro purchases due to Allocation Shifting at LPVED, with battery storage requirements.**

| Peak Shaving capability (MW) | Increase in Firm Hydro Purchases (MWh) | BESS energy Required (MWh) |
|---|---|---|
| 2018 | | |
| 0.5 | 1098 | 0.71 |
| 1.0 | 2224 | 2.03 |
| 2.0 | 4542 | 13.9 |
| 2019 | | |
| 0.5 | 1125 | 0.72 |
| 1.0 | 2279 | 2.16 |
| 2.0 | 4241 | 8.85 |

**Table 5. Peak shaving capability for LPVED in Summer Months.**

| Month | 0.5MW, 0.75MWh unit | 1.0MW, 2.20MWh unit | 2.0MW, 13.9MWh unit |
|---|---|---|---|
| May 2018 | 0.5MW | 1.0MW | 2.0MW |
| June 2018 | 0.38MW | 0.72MW | 1.85MW |
| July 2018 | 0.48MW | 0.69MW | 1.62MW |
| Aug. 2018 | 0.27MW | 0.43MW | 1.37MW |
| Sept. 2018 | 0.5MW | 1.0MW | 2.0MW |
| Oct. 2018 | 0.5MW | 1.0MW | 2.0MW |

### B. Peak Shaving to Reduce Demand Charge

In this study, battery energy ratings were chosen to meet the Table 4 requirements for allocation shifting. The chosen sizes are given in Table 5. Table 5 also shows the resulting peak shaving capability in months where the firm hydro allocation is not exceeded.

### C. LPVED Benefits and Costs

The projected benefits to LPVED from the selected BESS are shown in Table 6. Again, the 2.0 MW unit was not selected for further study due to its energy rating requirement. Table 6 shows that both the 0.5MW and 1.0MW units would be near the breakeven point at present, based on the projected savings and the annualized cost of the system for 2020.

**Table 6. Comparison for the Annual Benefits and Annualized Costs for the LPVED Peak Shaving Scenarios.**

| BESS Unit | Allocat. Shifting Benefit | Demand Reduction Benefit | Market Rate Optim. Benefit | Total Proj'ted Benefit | 2020 NaS system cost |
|---|---|---|---|---|---|
| 0.5MW | $22731 | $11641 | $1875 | $36247 | $23401 |
| 1.0MW | $46028 | $19698 | $3750 | $69476 | $65088 |

## VI. CASE 3: MASSENA ELECTRIC DEPARTMENT

The Massena Electric Department (MED) NYPA firm hydro allocation was 23.556 MW in 2018. In 2018, MED demand exceeded its allocation during all twelve months of the year. In 2019, demand exceeded the allocation in 11 of the 12 months.

### A. Allocation Shifting

Table 7 shows the allocation shifting results for MED. Both the energy shifted and BESS energy ratings are significantly higher than LPVED, primarily due to MED exceeding its allocation in all 12 months of 2018. Again, there is good consistency in the results between 2018 and 2019.

**Table 7. MED Allocation Shifting Results**

| Peak Shaving capability (MW) | Increase in Firm Hydro Purchases (MWh) | BESS energy Required (MWh) |
|---|---|---|
| **2018** | | |
| *0.5* | 2267 | 1.80 |
| *1.0* | 4609 | 4.50 |
| *2.0* | 9010 | 16.00 |
| **2019** | | |
| *0.5* | 2181 | 1.05 |
| *1.0* | 4435 | 4.08 |
| *2.0* | 8478 | 14.90 |

### B. Benefit/Cost Analysis

MED will not benefit from NYPA demand reduction in 2018 as the firm hydro allocation was exceeded in all 12 months of that year. However, NYPA also wheels the power directly to MED, and these wheeling charges have a demand based component.

Table 8 shows the projected benefits and costs for the MED 0.5MW and 1.0 MW installations. Table 8 shows that the projected annual savings exceed the NaS annualized cost projections in both cases.

**Table 8. Comparison for the Annual Benefits and Annualized Costs for the MED Peak Shaving Scenarios.**

| BESS Unit | Allocat. Shifting Benefit | Demand Reduction Benefit | Market Rate Optim. Benefit | Total Proj'ted Benefit | 2020 NaS system cost |
|---|---|---|---|---|---|
| **0.5MW** | $37196 | $14880 | $3750 | $55826 | $50829 |
| **1.0MW** | $75539 | $29760 | $7500 | $112800 | $109495 |

## VII. CONCLUSIONS

The project results show that battery energy storage systems can have a role in reducing overall costs to New York's municipal electric departments. The savings come from allocation shifting, peak demand reduction, and market arbitrage. Three case studies are presented in the paper. These case studies involve a range of municipal department sizes, relative to their firm hydro allocation. In these cases, the benefits were relatively higher when the firm hydro allocation was more often exceeded.

The report also notes that the impact of the wheeling charges will vary with the provider of that service.

## VIII. ACKNOWLEDGMENTS

This project would not have been possible without the close cooperation of: Marc Staves and Michael Dominie, Tupper Lake, Kimball Daby, Lake Placid and Andrew McMahon and Jeff Dobbins, Massena, and colleagues at their respective organizations. We acknowledge the efforts of Nancy Bernstein, ANCA, who saw the need for this project and brought the team together. We also acknowledge Darryl Jacobs on NYPA, for his insights on the rate structures for the municipal departments.